\documentclass[preprint,12pt]{elsarticle}
\newcounter{bla}

\begin{document}

\begin{frontmatter}

\title{Numerical Integration as an Initial Value Problem}
 \author{Daniel Gebremedhin}
 \ead{daniel1.gebremedhin@famu.edu}
 \author{Charles Weatherford}
 \ead{charles.weatherford@famu.edu}
\address{Physics Department, Florida A\&M University, Tallahassee, FL, USA.}
\date{\today}
\begin{abstract}
  Numerical integration (NI) packages commonly used in scientific research are
  limited to returning the value of a definite integral at the upper integration
  limit, also commonly referred to as numerical quadrature. These quadrature
  algorithms are typically of a fixed accuracy and have only limited ability to
  adapt to the application. In this article, we will present a highly adaptive
  algorithm that not only can efficiently compute definite integrals encountered
  in physical problems but also can be applied to other problems such as
  indefinite integrals, integral equations and linear and non--linear eigenvalue
  problems. More specifically, a finite element based algorithm is presented
  that numerically solves first order ordinary differential equations (ODE) by
  propagating the solution function from a given initial value (lower
  integration value). The algorithm incorporates powerful techniques including,
  adaptive step size choice of elements, local error checking and enforces
  continuity of both the integral and the integrand across consecutive
  elements.
\end{abstract}

\end{frontmatter}

\section{Introduction}
\label{sec:introduction}

Numerical integration(NI) is one of the most useful numerical tools that is
routinely utilized in all scientific and engineering applications. There are
general and specialized methods  that efficiently compute definite integrals
even for many pathological cases. Most of the integration algorithms are in
software packages and are available in almost all scientific
libraries. Ref.~\cite{davis2007methods} is an excellent book on the subject
containing a comprehensive reference to scientific articles and the accompanying
software. However, as efficient and advanced as those quadrature algorithms are,
and perhaps for this very reason, the techniques implemented in those algorithms
do not necessarily carry over to treat other related numerical problems.

Applications pertinent to NI are usually a lot more complex than simply a value
at a single point at the upper integration limit. A broader and versatile
advantage can be gained if integration is viewed from its general mathematical
perspective as a specific case of an initial value problem. In this paper, we
present an algorithm that solves ODEs that, with a slight modification, can be
used on many relevant applications, one of the most important of which being
NI. As far as numerical quadrature is concerned, the efficiency of the method
given here is comparable to the common NI packages available in scientific
libraries. Moreover, the present method is so simple to use that it can be
readily modified and applied to a broader spectrum of numerical applications as
long as they can be set up as initial value problems.

The algorithm in question has recently been applied to a second order ODE to
solve and treat a challenging system--namely, the soft Coulomb problem, where its
versatility and power is self--evident~\cite{PhysRevE.89.053319}. In this work,
we discuss in some detail how the same algorithm can be adapted to solve a
first-order ODE and thereby apply it, in a purely mathematical setting, to the
calculation of the numerical integral of a given function. Illustrative examples
that are best solved by the present method, and thereby highlight its important
features, will also be included.

\section{Description of Algorithm}
\label{sec:descr-algor}

The main intent of this paper is to present a finite element           algorithm that
numerically approximates a solution function $y$ to the following first order
ODE
\begin{equation}
  \label{eq:1}
  \frac{{\rm{d}}}{{\rm{d}}x} y(x) = f(x)
\end{equation}
\noindent with a given initial condition $y(x=a) = y_a$, where, both $y_a$ and
$a$ are assumed to be finite. Generally, $f(x)$ can be a simple function of $y$
such that the above equation may be a linear or non--linear eigenvalue
problem. $f$ can also be a kernel function of homogeneous or inhomogeneous
integral equation~\cite{kythe2002computational}. Particularly, if $f(x)$ is a
predetermined simple function and $y_a = 0$, then the above equation will be
equivalent to a single integral given by
\begin{equation}
  \label{eq:2}
  y(x) = \int_a^x f(t) \, {\rm{d}}t.
\end{equation}

Begin the numerical solution to eq.~(\ref{eq:1}) by breaking the x--axis
into finite elements and mapping into a local variable $\tau$ with domain $-1
\leq \tau \leq 1$ defined by a linear transformation given below.
\begin{equation}
  \label{eq:3}
  x = x_i + q_i(\tau + 1), \qquad x_i \le x \le x_{i + 1}
\end{equation}
Here, $i = 1, 2, \ldots, i_{\rm{max}}$ labels the elements with $x_1 = a$, while
$q_i = (x_{i + 1} - x_i)/2$ is half the size of the element. In terms of the
local variable eq.~(\ref{eq:1}) can be re--written as
\begin{equation}
  \label{eq:4}
  \frac{{\rm{d}}}{{\rm{d}}\tau} \bar{y}(\tau) = q \bar{f}(\tau).
\end{equation}
The over--bar indicates the appropriate change in functional form while the
element index $i$ is dropped for simplification of notation. At this point, we
will expand $\bar{y}(\tau)$ in a polynomial basis set as follows.
\begin{equation}
  \label{eq:5}
  \bar{y}(\tau) = \sum_{\mu = 0}^{M-1} u_\mu(\tau)B_\mu + s_0(\tau) q
  \bar{f}(-1) + \bar{y}(-1)
\end{equation}
\noindent Notice that $\bar{f}(-1) \equiv f(x_i)$ and similarly for $\bar{y}$. The
main import of the above expansion is that it allows us to enforce continuity of
both $f$ and $y$ across the boundary of two consecutive elements. This is
because the basis functions $u$ and their derivatives (denoted by $s$)
identically vanish at $-1$. These functions $u$ and $s$ are defined in terms of
Legendre polynomials of the first kind $P$ \cite{citeulike:1816367} as
\begin{equation}
  \label{eq:6}
  s_{\mu}(\tau) = \int_{-1}^{\tau} P_{\mu}(t) \, \rm{d}t, \qquad u_{\mu}(\tau) =
  \int_{-1}^{\tau} s_{\mu}(t) \, \rm{d}t.
\end{equation}
\noindent
They are discussed in~\cite{PhysRevE.89.053319} in more detail. These
polynomials satisfy the following recurrence identities presented here for the
first time. It is surprising to see how these relations remain three--term, with
no surface values, despite the fact that the polynomials $s$ and $u$ are
sequentially primitives of  Legendre polynomials.
\begin{eqnarray}
  \label{eq:7}
  s_0(\tau) &=& \tau + 1, \qquad s_1(\tau) = \frac{1}{2}(\tau^2-1) \nonumber \\
  (\mu + 1)s_\mu(\tau) &=& (2\mu - 1)\tau s_{\mu - 1}(\tau) - (\mu - 2)s_{\mu -
                           2}(\tau) \qquad \mu \ge 2 \\
  \label{eq:8}
  u_0(\tau) &=& \frac{1}{2}(\tau+1)^2, \qquad u_1(\tau) =
                \frac{1}{6}(\tau+1)^2(\tau-2) \nonumber \\
  (\mu + 2)u_\mu(\tau) &=& (2\mu - 1)\tau u_{\mu - 1}(\tau) - (\mu - 3)u_{\mu -
                           2}(\tau) \qquad \mu \ge 2
\end{eqnarray}
\noindent The above two relations allow us to employ the Clenshaw recurrence
formula~\cite{Press:1992:NRF:141273} which is known to facilitate effective
numerical evaluation of relevant summations as the one included in
eq.~(\ref{eq:5}).

Substituting the expansion given in eq.~(\ref{eq:5}) into eq.~(\ref{eq:4}),
evaluating the resulting equation at Gauss--Legendre (GL) abscissas and after
some rearrangement, we get the following set of simultaneous equations of size
$M$,
\begin{equation}
  \label{eq:9}
  \sum_{\mu = 0}^{M-1} s_\mu(\tau_{\nu})B_\mu = q \left[ \bar{f}(\tau_{\nu}) -
    \bar{f}(-1) \right]
\end{equation}
\noindent where, $\tau_{\nu}$ is a root of an $M^{\rm{th}}$ order Legendre
polynomial. This technique, known as collocation method, is an alternate way of
constructing a linear system of equations compared to the more familiar
projection integrals. $s_\mu(\tau_{\nu})$ are now elements of a square matrix
which, for a given size $M$, is constant and hence, merely needs to be
constructed only once, LU decomposed and stored for all times. Thus solving for
the unknown coefficients $B$ only involves back substitution. Construction of
the right--hand side column, on the other hand, requires an evaluation of the
function $f$ at $M + 1$ points including at the beginning of the element. In the
examples that follow we will denote the total number of function evaluations as
$N$.

After the coefficients $B$ are calculated this way, the value of the integral at
the end point of the element $i$ can be obtained by evaluating eq.~(\ref{eq:5})
at $\tau=+1$. The result is
simply
\begin{equation}
  \label{eq:10}
  \bar{y}(+1) = 2B_0 - \frac{2}{3}B_1 + 2q\bar{f}(-1) + \bar{y}(-1).
\end{equation}
\noindent
One of the advantages of representing the solution in the form given in
eq.~(\ref{eq:5}) is that it allows us to evaluate the value of the integral at
any continuous point inside the element. In this sense, $y$ is simply a function
whose domain extends into all of the solved elements as the propagation
proceeds. Hence, from a numerical perspective, this algorithm is best implemented
by way of object--oriented programming in order to incorporate and preserve
all the necessary quantities of all the elements in a hierarchy of derived
variables. The integrand function $y(x)$, for instance, can be saved (as an
object) for later use by including, among other quantities, the grid containing
the coordinates of the steps $\left\{ x_i \right\}_{i = 1}^{i_{\rm{max}}}$. Then
locating the index to which a given point $x$ belongs, is an interpolation
exercise for which efficient codes already exist. See, for instance, subroutine
\emph{hunt} and its accompanying notes in ref.~\cite{Press:1992:NRF:141273}. For
a sorted array, as is our case, this search generally takes about $\log{}_2
i_{\rm{max}}$ tries while further points within close proximity can be located
rather easily.

An even more important advantage of this form of the solution is its suitability
for estimation of the size of the next element via the method described
in~\cite{PhysRevE.89.053319}. This adaptive step size choice relies on the
knowledge of the derivatives of $y$ up to $4^{\rm{th}}$ order at the end of a
previous element. Those derivatives can be directly computed from
eq.~(\ref{eq:5}). Instead of solving two quadratic equations as we did in the
last paper, however, we will solve one cubic equation, which is more suited for
the present purpose. The user is required to guess only the size of the very
first element. Overestimation of the size of an element may cause an increase in
the number of function evaluations, but there will be no compromise in the
accuracy of the solution as it will be made clear in a moment.

One other attractive feature is that we can precisely measure the error of the
calculated solution directly from eq.~(\ref{eq:4}). Specifically, at the end of
a solved element $i$, the error between the exact $f(x_{i + 1})$ and
the calculated
\begin{eqnarray}
  \label{eq:11}
  \frac{\rm{d}}{{\rm{d}}x} y(x)\biggr|_{x=x_{i+1}} = \frac{1}{q_i}
  \frac{\rm{d}}{{\rm{d}} \tau} \bar{y}(\tau) \biggr|_{\tau=+1}
  = \frac{2B_0}{q_i} + f(x_i)
\end{eqnarray}
\noindent can be obtained. Notice that $f(x_{i + 1})$ will be needed in the
beginning point of the next element. If the resulting error is not satisfactory,
a bisection step will be taken to reduce the size of the element by
moving the upper limit $x_{i + 1}$ closer to $x_i$ and re--solving. The main
purpose of the adaptive step size choice is then to estimate \emph{a priory} an
optimum step size, thereby reducing the number of bisections and/or the number
of function evaluations necessary. Since this error is that of the integrand
$f$, and not of the integral $y$, it is not necessary to demand this error be as
small as the machine precision. The reason is, higher derivatives of $y$
computed from eq.~(\ref{eq:5}) will successively deteriorate in accuracy as the
order increases. In a $16$ digit calculation a relative error of $10^{-4} -
10^{-7}$ in $f(x_{i + 1})$ often suffices, depending on how smooth
the integrand is, for calculating the solution function $y$ correctly to within
$\mathcal{O}(-14)$. Other quadrature methods do not directly take the integrand
into account in their error estimation. They often use a formula, an outcome
usually of a non--trivial analytic derivation, that estimates the
upper bound of a residual term for a specific order~\cite{Ehrich1999}.

In the present algorithm, the number of basis functions is kept constant. Other
methods such as Clenshaw--Curtis quadrature, take advantage of a convenient
property of the roots of Chebyshev polynomials that allows for conserving
preceding function evaluations whenever the order of expansion is doubled. But
this recursive doubling of basis set size is not necessarily effective as far as
improving the accuracy of the integration is concerned. The reason is, for a
given working precision, there is usually only a limited range of an optimum
number of basis functions, say $10-16$, whose half or double is either too small
or excessive. Rather, a more effective way is to estimate an adaptive step size,
fix an optimum order of expansion, and reduce the size of the element whenever
necessary - which is what is done here.

Implementation of the algorithm starts by fixing the size of the first step
$q_1$ and the size of the linear system $M$. Also, define $f(a) = f(x_1)$. The
rest of the procedure is sketched below.
\begin{enumerate}
\item Construct the right--hand side of eq.~(\ref{eq:9}) by evaluating the
  function  $\bar{f}(\tau_\nu)$ at the GL nodes.
\item Calculate the solution coefficients $B$ by back substitution.
\item Calculate the error between the two sides of eq.~(\ref{eq:2}) at the end
  point using eq.~(\ref{eq:11}) for the left--hand side and by evaluating
  $f(x_{i+1})$ directly. Retain the value of $f(x_{i+1})$ which will be needed
  at the beginning of the next element.
\item If the error is not satisfactory, reduce $q_i \rightarrow q_i/2$ and go
  to step~(1).
\item Otherwise calculate $y$ as per eq.~(\ref{eq:10}) and its higher
  derivatives  at $\tau = 1$ using  eq.~(\ref{eq:5}) and estimate the size of
  the next step $q_{i+1}$ using the method presented
  in~\cite{PhysRevE.89.053319}.
\item For finite domain systems with an upper limit $b$, make sure the resulting
  step size does not stride beyond $b$. i.e., take $q_{i+1} \leftarrow
  min(q_{i+1}, (b - x_{i+1})/2)$. For open integrals on the other hand, keep
  propagating until the value of the integral $y$ converges.
\end{enumerate}

Finally, we may sometimes prefer not to evaluate the function $f(a)$ at the
lower limit of the integration. A good example is an integrand containing
singular terms at the origin. In those instances, only for the very first
element, the main expansion given in eq.~(\ref{eq:5}) can be altered to be in
terms of the $s$ polynomials as follows.
\begin{equation}
  \label{eq:12}
  \bar{y}(\tau) = \sum_{\mu = 0}^{M-1} s_\mu(\tau)B_\mu + \bar{y}(-1)
\end{equation}
\noindent The rest of the propagation can then
resume normally and the other modifications that follow can be worked out
straightforwardly.

\section{Numerical Examples}
\label{sec:numerical-examples}

We will now consider illustrative examples that demonstrate the typical behavior
of the numerical algorithm discussed above. All of the examples chosen are
difficult to calculate with an accuracy close to working precision, without
breaking the range of integration into smaller intervals. Whenever applicable,
comparisons will be made with DQAG~\cite{Quadpackref}, which is one of the
integration subroutines compiled in QUADPACK
\footnote{http://www.netlib.org/quadpack}. DQAG is an adaptive method that keeps
bisecting the element with highest error estimate until the value of the overall
integration is achieved to within the desired error.

In all of the examples that follow, the error at the end of the elements has
been determined by
\begin{eqnarray}
  \label{eq:13}
  \left| \frac{2B_0}{q_i} + f(x_i) - f(x_{i+1}) \right| \le \left| f(x_{i+1})
  \right| \delta_{\rm{rel}} + \delta_{\rm{abs}}
\end{eqnarray}
\noindent where, $\delta_{\rm{abs}}$ and $\delta_{\rm{rel}}$ denote the absolute and
relative errors respectively. The values of $\delta_{\rm{rel}} \approx 2.22
\times 10^{-4}$ and $\delta_{\rm{abs}} \approx 2.22 \times 10^{-19}$ will be used
unless otherwise specified. The size of the first step is fixed to be 0.5--that is
$q_1 = 0.25$ and the number of basis functions used is $M = 13$.

DQAG, on the other hand, has been run with absolute and relative errors
of $\delta_{\rm{rel}} = \delta_{\rm{abs}} \approx 1.11 \times 10^{-13}$ and
the lowest order ($key = 1$) has been used in order to keep the number of
function evaluations at a minimum so as to provide a fair comparison. Our program is written
in modern C++ and run on a late 2013, 2.8 GHz Intel Core i7 MacBook Pro laptop
computer with Apple LLVM 8.1 compiler.

\subsection{Closed Integral}
\label{sec:closed-integral}

We consider the $15$ closed integrals studied in
~\cite{doi:10.1080/10586458.2005.10128931}. These integrals exhibit different
properties and were chosen by Bailey \emph{et al.} to test three
other methods of NI. Herein we attempt all of them and report the
output. Table~\ref{tab:table1} shows the results for numerical values of those
integrals calculated using the present method as well as DQAG.
\begin{table}[htp]
  \caption{Results for $15$ closed integrals taken
    from~\cite{doi:10.1080/10586458.2005.10128931}. The values of the integral
    $y(b)$ from the present method are shown in the second
    column. $\delta_{\rm{rel}}$ and $N$ of both the present method and DQAG are
    also shown.}
  \begin{center}\footnotesize
    \renewcommand{\arraystretch}{1.3}
    \begin{tabular}{|c|c|l|l|l|l|}\hline
      \textrm{label}  &  \textrm{$y(b)$}  &  $\delta_{\rm{rel}}$  &
      $\delta_{\rm{rel}}$ (DQAG) &\multicolumn{1}{|c|}{$N$}
      &  $N$ (DQAG) \\
      \hline
      $1$   &  0.250 000 000 000 000   &  1.110[-16]  &  0.000       &  29    &
                                                                                15 \\
      $2$   &  0.210 657 251 225 807   &  0.000       &  1.318[-16]  &  29    &
                                                                                45 \\
      $3$   &  1.905 238 690 482 68    &  0.000       &  0.000       &  191   &
                                                                                15 \\
      $4$   &  0.514 041 895 890 071   &  0.000       &  0.000       &  29    &
                                                                                45 \\
      $5$   &  -0.444 444 444 444 445  &  8.743[-16]  &  3.747[-16]  &  871   &
                                                                                885 \\
      $6$   &  0.785 398 163 397 448   &  1.414[-16]  &  2.827[-16]  &  974   &
                                                                                795 \\
      $7$   &  1.198 140 227 142 81    &  6.337[-9]   &  6.159[-9]   &  2129  &
                                                                                1725 \\
      $8$   &  1.999 999 999 999 98    &  8.438[-15]  &  4.441[-16]  &  922   &
                                                                                1545 \\
      $9$   &  -1.088 793 045 151 79   &  9.993[-15]  &  2.855[-15]  &  1243  &
                                                                                1335 \\
      $10$  &  2.221 441 454 672 65    &  6.485[-9]   &  2.907[-9]   &  2032  &
                                                                                1725 \\
      $11$  &  1.570 796 326 794 90    &  0.000       &  0.000       &  29    &  105 \\
      $12$  &  1.772 453 840 168 93    &  6.057[-9]   &  5.929[-9]   &  2439  &
                                                                                1455 \\
      $13$  &  1.253 314 137 315 62    &  9.514[-14]  &  0.000       &  96    &  255 \\
      $14$  &  0.500 000 000 000 001   &  1.110[-15]  &  0.000       &  231   &  375 \\
      $15$  &  1.570 796 326 794 90    &  1.414[-16]  &  2.218[-11]  &  1523  &
                                                                                210 \\
      \hline
    \end{tabular}
  \end{center}
  \label{tab:table1}
\end{table}
$\delta_{\rm{rel}}$ with respect to the exact values, and a number of function
evaluations $N$ are also indicated. $\delta_{\rm{rel}}$ and $N$ of the output
from DQAG are also shown for comparison. The two methods are essentially,
qualitatively similar in performance. Notice that the labels of the integrals
are taken from the article ~\cite{doi:10.1080/10586458.2005.10128931}.

\begin{figure}
  \includegraphics[scale=0.7]{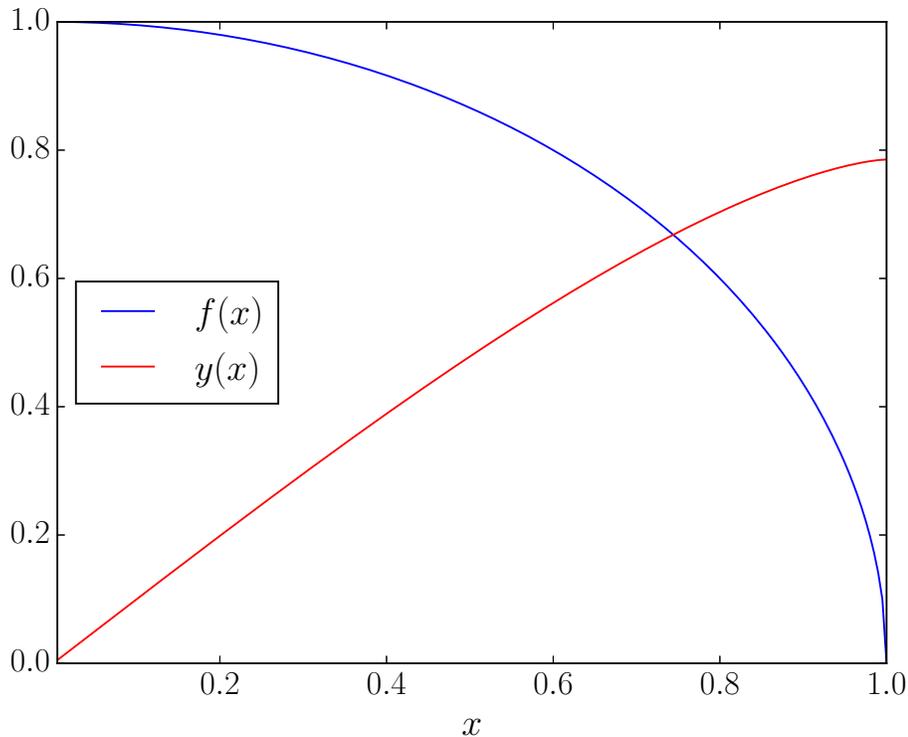}
  \caption{(Color online) Plots of the functions given in Eqs.~(\ref{eq:14})
    and~(\ref{eq:15}). $f$ varies very rapidly towards the upper integration
    limit causing numerical difficulty.}
  \label{fig:1}
\end{figure}

\begin{figure}
  \includegraphics[scale=0.8]{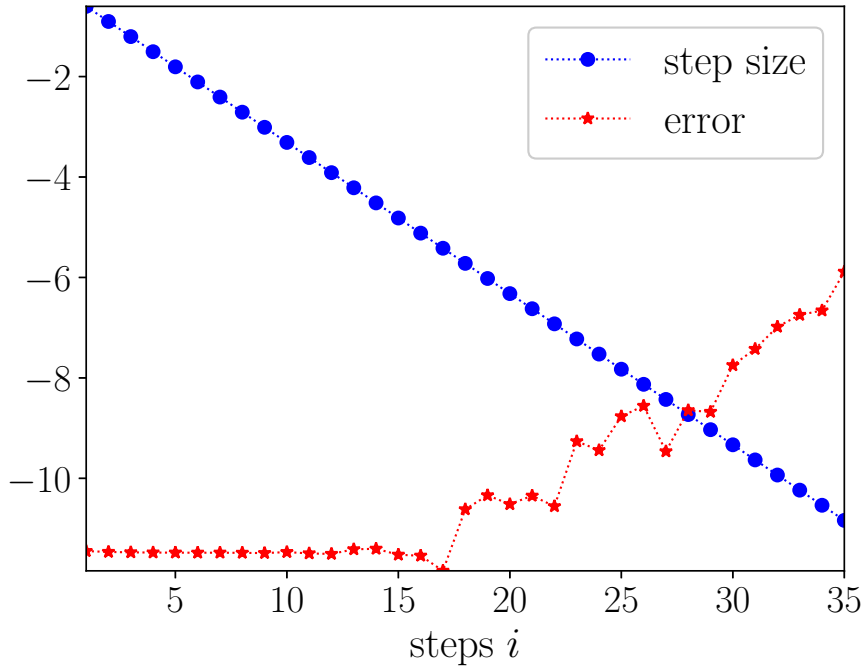}
  \caption{(Color online) Step size of the finite elements (blue .) and
    $\delta_{\rm{rel}}$ of the integral at the end of the elements (red *) are
    shown. The algorithm spends more time and picks up more errors towards the
    upper integration limit. Vertical axis values are $\log{}_{10}$ of the
    actual.}
  \label{fig:2}
\end{figure}

Lets take a closer look at one of the integrals (index 6) in order to
demonstrate the behavior of our algorithm. Example 6 has an integration limits
$[0, 1]$ and an integrand $f$ given below.
\begin{equation}
  \label{eq:14}
  f(x) = \sqrt{1 - x^2}, \qquad 0 \le x \le 1
\end{equation}
\noindent
The integral of the above function is
\begin{equation}
  \label{eq:15}
  y(x) = \frac{1}{2} \left[ x \sqrt{1 - x^2} + \arcsin(x) \right].
\end{equation}
\noindent with $y(1) = \pi/4$. Both functions $f$ and $y$ are plotted in
Fig.~\ref{fig:1}.
\noindent The calculated value of the integral is $y(1) = 0.785 398 163 397 448$
which has a $\delta_{\rm{rel}}$ of $\approx 1.414\times 10^{-16}$ compared to
the exact value. It took $N = 974$ function evaluations and $i_{\rm{max}} = 35$
steps to propagate the integral $y$ from $0$ to $1$, where the size of the first
step was set to $0.5$ as mentioned above. As can be seen from
Fig.~(\ref{fig:2}), the step sizes chosen by our algorithm kept decreasing to as
low as $2.22\times 10^{-11}$, which is consistent with the singularity of the
derivative of the integrand $f$ at the upper integration limit.
Fig.~(\ref{fig:2}) also shows relative errors of the integral $y$ at the end of
all the elements in comparison with the exact values from eq.~(\ref{eq:15}).

Clearly, the maximum error occurs at the last step near $1$, which is what
motivated our choice of this example. Rapidly changing functions, such as those
with integrable singularity, are generally troublesome to the algorithm because
the step size may not be small enough and/or the order of polynomials high
enough to accommodate portions of the integrand with (nearly) vertical shape.

In comparison, DQAG takes $795$ function evaluations and calculates the integral
with a $\delta_{\rm{rel}}$ of $2.827\times 10^{-16}$.

\subsection{Nonlinear Problem}
\label{sec:nonlinear-problem}

Bender \emph{et al.} have studied the following interesting nonlinear eigenvalue
problem for which the numerical solution can be very
challenging~\cite{1751-8121-47-23-235204}.
\begin{equation}
  \label{eq:16}
  \frac{{\rm{d}}}{{\rm{d}}x} y(x) = \cos \left[\pi xy(x) \right], \qquad x \ge 0
\end{equation}
\noindent Since this is a nonlinear equation, it has to be solved in an
iterative fashion as:
\begin{equation}
  \label{eq:17}
  \frac{{\rm{d}}}{{\rm{d}}x} y_{\sigma + 1}(x) = f_{\sigma}(x), \qquad \sigma =
  1, 2, \ldots
\end{equation}
\noindent where $f_\sigma(x) = \cos \left[\pi xy_\sigma(x) \right]$ and $\sigma$
labels the levels of the iteration. Notice that at any stage of the iteration
only the values of $f$ at the GL nodes are required. The first iterate of array
$f_1$ is seeded from the solution vector $B$ of the final result in the
preceding element. This is not only convenient but also an excellent
approximation because the adaptive step choice implemented here is based on the
assumption that the solution function between two consecutive elements remains
constant up to a fourth order Taylor series expansion. The first element has
been started by setting all the elements of the solution vector $B$ to
unity. Apparently, only for this problem, the first two steps of the procedure
given in Section~\ref{sec:descr-algor}, must be repeated until the iteration in
eq.~(\ref{eq:17}) converges.

\begin{figure}
  \includegraphics[scale=0.7]{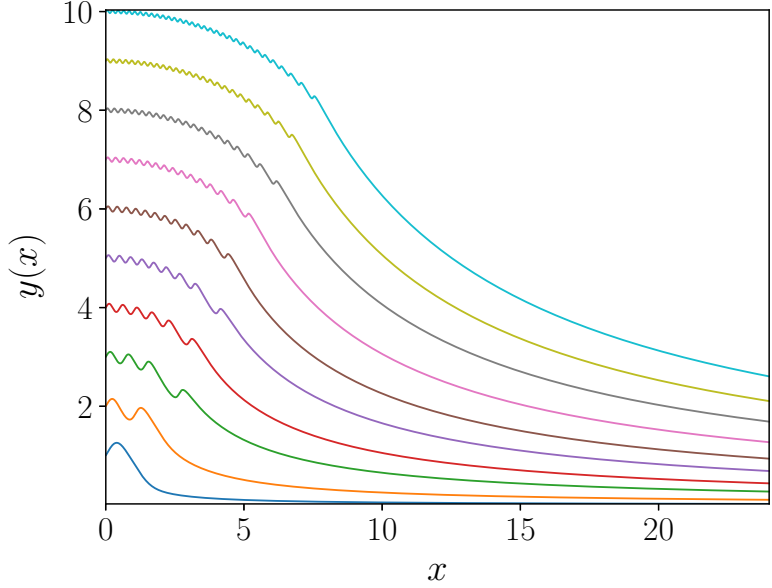}
  \caption{(Color online) Calculated solution functions to eq.~(\ref{eq:16}) for
    $0 \le x \le 24$ with initial condition $y_a = 1, 2, \ldots , 10$ are
    shown. The functions are oscillatory near the origin before they evolve
    into an asymptotic discrete bundle.}
  \label{fig:3}
\end{figure}

The results of our calculation for $y(x), x\in [0, 24]$, are plotted in
fig.~\ref{fig:3} for initial values at the origin $y_a = n, n = 1, 2, \ldots ,
10$.
Similar plots have been reported in ref.~\cite{1751-8121-47-23-235204} which are
a result of point--wise convergent calculations. Our solution function, on the
other hand, has been propagated from the origin outward. The very large number
of steps the solution required for such a modest distance of $x = 24$ is quite
remarkable. Table~\ref{tab:table1} summarizes the output of our program.
\begin{table}[htp]
  \caption{Results for $y(24)$ to eq.~(\ref{eq:16}) for the respective initial
    conditions $y(0)$. The rest of the columns are: total number of steps,
    average step size, average number of functions evaluations per element and
    time elapsed.}
  \begin{center}\footnotesize
    \renewcommand{\arraystretch}{1.3}
    \begin{tabular}{|c|c|c|c|l|c|}\hline
      \textrm{$y(0)$}  &  \textrm{$y(24)$}  & \textrm{$i_{\rm{max}}$} &
      \textrm{$2 q_{\rm{ave}} \, (\times 10^{-5})$}  &
      \multicolumn{1}{|c|}{\textrm{$N_{\rm{ave}}$}} &  \textrm{$t \, [sec]$} \\
      \hline
      $1$   &  0.020 844 865 419 015 3  &  1 633 376  &  1.469   &  254.338  &  63   \\
      $2$   &  0.104 224 327 270 128    &  1 701 378  &  1.411   &  253.336  &  66   \\
      $3$   &  0.270 983 253 633 302    &  1 908 989  &  1.257   &  251.829  &  73   \\
      $4$   &  0.437 742 187 280 245    &  2 068 413  &  1.160   &  250.609  &  79   \\
      $5$   &  0.687 880 611 222 152    &  2 397 070  &  1.001   &  249.380  &  91   \\
      $6$   &  0.938 019 076 811 230    &  2 636 633  &  0.9103  &  248.157  &  100  \\
      $7$   &  1.271 537 122 002 93     &  3 041 709  &  0.7890  &  247.071  &  116  \\
      $8$   &  1.688 434 875 810 57     &  3 586 933  &  0.6691  &  245.951  &  135  \\
      $9$   &  2.105 332 915 403 23     &  4 021 621  &  0.5968  &  244.952  &  151  \\
      $10$  &  2.605 611 041 676 66     &  4 626 563  &  0.5187  &  244.110  & 174 \\
      \hline
\end{tabular}
\end{center}
\label{tab:table2}
\end{table}
It took millions of steps, with an average step size as low as $\approx
5.187\times 10^{-6}$ and an average number of function evaluations up to
$N_{\rm{ave}} \approx 254.3$ per element. Notice that this $N_{\rm{ave}}$
includes all the iterations, which, along with the brevity of the elapsed times
shown, indicates the efficiency of our implementation.

We have also included the final numerical values at $y(24)$ for reference. For
this example, the required error has been set lower as $\delta_{\rm{rel}} =
3.0\times 10^{-9}$, for an obvious reason. In order to check the validity our
results, we have re--run the program by still lowering the magnitude of the
error to $\delta_{\rm{rel}} = 3.0\times 10^{-10}$. With the lowered error, the
indicated values of $y(24)$ vary by an amount no larger than $\sim 1.86\times
10^{-13}$. The computational times in the last column of the table also
increased to as high as $470$ seconds. Maintaining this much accuracy after
millions of steps, and in an iterative calculation, shows how robust our method
is. This non--linear eigenvalue problem is computationally challenging indeed.

\subsection{Double--Range Integrals}
\label{sec:double-range-integr}

In this example, we consider a double integral for which one of the limits
of the inner integral is identical to the variable of the outer integral. These
types of integrals are common in studies of many particle dynamical systems
involving Green's functions or double range addition theorems such as the
Laplace expansion. Particularly, we will look at the following integral, which is
the most prominent radial integral that is encountered in calculations that
involve exponential type orbitals or Geminals in a spherical coordinate
system. It stems from the addition theorem for ${r_{12}}^{n} e^{-\alpha r_{12}}$
given in~\cite{QUA:QUA24319} and the resulting integrals are still topics of
interest in recent research~\cite{Jiao2015140, Rico2012, QUA:QUA21002}. Let us
define the integral as
\begin{equation}
  \label{eq:18}
  \hspace*{-0.5cm} I_{\lambda_1 \lambda_2}^{\mu_1 \mu_2}(\alpha_1, \beta_1,
  \alpha_2, \beta_2) = \int_0^\infty {\rm{d}} y \, e^{-\alpha_1 y} y^{{\mu}_1}
  \hat{i}_{\lambda_1}(\beta_1 y) \int_y^\infty {\rm{d}} x \, e^{-\alpha_2 x}
  x^{{\mu}_2} \hat{k}_{\lambda_2}(\beta_2 x)
\end{equation}
\noindent where $\hat{i}$ and $\hat{k}$ are spherical modified Bessel functions
of the first and second kind respectively~\cite{citeulike:1816367}. The
screening parameters $\alpha_1, \beta_1, \alpha_2, \beta_2$ are positive real
numbers while all of the indices $\mu_1,\lambda_1, \mu_2, \lambda_2$ are
integers. The correct composition of these parameters is such that both the
inner and outer integrals remain finite as is the case in physical
applications.

In order to propagate from the origin, both lower integration limits need to be
set to zero. This can be attained by switching the order of the two integrals
using the following identity which maintains identical $x$--$y$ region of
integration~\cite{hassani2000mathematical}.
\begin{equation}
  \label{eq:19}
  \int_0^\infty {\rm{d}} y \, f(y) \int_y^\infty {\rm{d}} x \, g(x) \equiv
  \int_0^\infty {\rm{d}} x \, g(x) \int_0^x {\rm{d}} y \, f(y)
\end{equation}
\noindent Hence, eq.~(\ref{eq:18}) can be written as
\begin{equation}
  \label{eq:20}
  \hspace*{-0.5cm} I_{\lambda_1 \lambda_2}^{\mu_1 \mu_2}(\alpha_1, \beta_1,
  \alpha_2, \beta_2) = \int_0^\infty {\rm{d}} x \, e^{-\alpha_2 x} x^{{\mu}_2}
  \hat{k}_{\lambda_2}(\beta_2 x) J_{\lambda_1}^{\mu_1}(\alpha_1, \beta_1; x)
\end{equation}
\noindent where, $J$ now represents the inner integral given below which needs
to be propagated only once from the origin until the integral converges.
\begin{equation}
  \label{eq:21}
J_{\lambda_1}^{\mu_1}(\alpha_1, \beta_1; x) = \int_0^x {\rm{d}} y \,
e^{-\alpha_1 y} y^{{\mu}_1} \hat{i}_{\lambda_1}(\beta_1 y), \qquad 0 \le x \le
x_{i_{\rm{max}}+1}
\end{equation}
\noindent Here $x_{i_{\rm{max}}+1}$ signifies the end point of the last element
where the result of the above integral converged to its value at
infinity. Hence, beyond this point $J$ is considered to be constant function,
i.e., $J(x) = J(x_{i_{\rm{max}}+1})$ for $x>x_{i_{\rm{max}}+1}$. Once the above
integral is done and all the relevant parameters stored, it can be evaluated at
any desired point $x>0$, which is a significant computational gain since the
double integral $I$ given in eq.~(\ref{eq:18}) has essentially been reduced to
two simple integrals. This demonstrates one of the main advantages contained in
the present algorithm.

Table~\ref{tab:table2} shows a sample calculation for the integral $I$ in
eq.~(\ref{eq:20}) for $\beta_1, \beta_2 \in \{0.5, 1.0, 2.0 \}$.
\begin{table}[htp]
  \caption{Exact and calculated values of the integral $I$ are shown for
    $\beta_1$ \& $\beta_2 \in \{0.5, 1.0, 2.0 \}$. Columns 3 and 4 show the
    number of function evaluations $N$ taken in the integrals Eq.~(\ref{eq:20})
    and Eq.~(\ref{eq:19}) respectively. Numbers in square bracket signify powers
    of ten.}
\begin{center}\footnotesize
\renewcommand{\arraystretch}{1.3}
\begin{tabular}{|c|c|c|c|l|l|}\hline
  \textrm{$\beta_1$}  &  \textrm{$\beta_2$}  &  \textrm{N of $J$}
  &  \textrm{N of $I$}  &  \multicolumn{1}{|c|}{I -- calculated}
  &  \multicolumn{1}{|c|}{I -- exact} \\
  \hline
      $0.5$  &  0.5  & 219  &  259  &  1.627 473 168 386 65 [27]  &  1.627 473 168 386 653 87 [27] \\
      $0.5$  &  1.0  & 219  &  218  &  2.559 085 779 949 79 [22]  &  2.559 085 779 949 794 01 [22] \\
      $0.5$  &  2.0  & 219  &  231  &  3.103 777 873 917 21 [17]  &  3.103 777 873 917 210 86 [17] \\
      $1.0$  &  0.5  & 232  &  259  &  2.946 389 365 576 82 [23]  &  2.946 389 365 576 741 23 [23] \\
      $1.0$  &  1.0  & 232  &  245  &  6.062 810 005 197 87 [18]  &  6.062 810 005 197 874 73 [18] \\
      $1.0$  &  2.0  & 232  &  204  &  9.533 337 428 978 80 [13]  &  9.533 337 428 978 808 27 [13] \\
      $2.0$  &  0.5  & 245  &  259  &  4.342 544 722 241 59 [19]  &  4.342 544 722 241 718 83 [19] \\
      $2.0$  &  1.0  & 245  &  245  &  1.097 615 571 907 40 [15]  &  1.097 615 571 907 438 80 [15] \\
      $2.0$  &  2.0  & 245  &  231  &  2.258 572 729 378 12 [10]  &  2.258 572 729 378 146 95 [10] \\  \hline
\end{tabular}
\end{center}
\label{tab:table2}
\end{table}
The rest of the  parameters are set as $\lambda_1 = -11, \mu_1 = 12, \lambda_2 =
-13, \mu_2 = 14, \alpha_1 = 2\beta_1, \alpha_2 = 2\beta_2$. The Bessel functions
$\hat{i}$ and $\hat{k}$ have been computed using the subroutines in GNU
Scientific Library (GSL) \footnote{http://www.gnu.org/software/gsl/}. For the
first element, eq.~(\ref{eq:12}) has been used in order to avoid evaluation at
the origin. We also have calculated the exact values of the integrals using
\emph{Mathematica} \cite{Mathematica7}, the first 18 digits of which are
displayed in the last column. Comparison with the the present calculated results
reveals that the integral $I$ was done accurately. The table also further shows
the total number of function evaluations $N$ for the inner and outer integrals.

\section{Conclusion}
\label{sec:conclusion}

There are many physical applications that can be modeled as initial value
problems and the methods available to compute them are equally diverse. No
single algorithm is known to address all of them at once, and hence, researchers
usually digress from their main area of interest so as to familiarize themselves
with many of the computational options. The present algorithm is by no means
capable of solving all initial value problems, but it comes pragmatically close,
especially for most physical applications. This is possible because it produces
solutions on finite elements whose size is chosen to locally, checks the
validity of the solution, and communicates the solution function and its first
derivative to the next element, maintaining continiuity of both. More
importantly, it can be applied on any ODE because it is easy to implement and
modify, especially with the proper use of object--oriented programming
techniques. The basis functions $u$ and $s$ discussed above are based on
Legendre polynomials. In the future, we will use other classic orthogonal
polynomials such as Chebyshev or Jacobi to determine if further advantages can
be gained.

\section{Acknowledgement}
\label{sec:acknowledgement}
DHG and CAW were partially supported by the Department of Energy, National
Nuclear Security Administration, under Award Number(s) DE-NA0002630. CAW was
also supported in part by the Defense Threat Reduction Agency.
\section{REFERENCES}
\label{sec:references}
\bibliographystyle{elsarticle-num}
\bibliography{SPTaQuadrature}
\end{document}